\documentclass{article}
\usepackage{amsmath}
\usepackage{amsfonts}   
\usepackage{amssymb}

\newcommand{\mlabel}[1]{\label{#1}
}

\newcommand{\seq}{\begin{equation}}                 
\newcommand{\eeq}[1]{\label{#1}\end{equation}
    }

\newcommand{\pf}{ \par \vspace{1ex} \noindent {\sc Proof} \hspace{2mm}}
\newcommand{\epf}{$ \quad \Box$ \par \vspace{1ex}}

\newtheorem{Theorem}{Theorem}[section]
\newcommand{\sthm}{\begin{Theorem}}         
\newcommand{\ethm}{\end{Theorem}}           

\newtheorem{Corollary}[Theorem]{Corollary}
\newcommand{\scor}{\begin{Corollary}}       
\newcommand{\ecor}{\end{Corollary}}         

\newtheorem{Lemma}[Theorem]{Lemma}
\newcommand{\slm}{\begin{Lemma}}            
\newcommand{\elm}{\end{Lemma}}              

\newtheorem{Remark}[Theorem]{Remark}
\newcommand{\srmark}{\begin{Remark}\rm}        
\newcommand{\ermark}{\end{Remark}}             
\newtheorem{Example}[Theorem]{\sc Example}
\newcommand{\sex}{\begin{Example}\rm}        
\newcommand{\eex}{\end{Example}}             

\newcommand{\seql}{\begin{eqnarray*}}       
\newcommand{\eeql}{\end{eqnarray*}}

\newcommand{\smlist}[1]{\begin{list}           
                      {(#1{zzcount})}{\usecounter{zzcount}}}
\newcommand{\elist}{\end{list}}

\newcommand{\for}{\quad \mbox{\textrm{for}} \quad}

\newcommand{\e}{\varepsilon}

\newcommand{\G}{\Gamma}
\renewcommand{\l}{\lambda}

\newcommand{\s}{\sigma}

\renewcommand{\t}{\tau}
\renewcommand{\O}{\Omega}

\newcommand{\vect}[1]{\mathbf{#1}}
\providecommand{\abs}[1]{\lvert#1\rvert}
\providecommand{\norm}[1]{\lVert#1\rVert}

\begin{document}
\pagenumbering{arabic}
\title{\bf{\small{\uppercase{Radial Convex Solutions of Boundary Value Problems for Systems of Monge-Amp\`ere equations}}}}

\author{Haiyan Wang\\\\Division of Mathematical and Natural Sciences \\ Arizona State University\\Phoenix, AZ 85069-7100, U.S.A.\\E-mail: wangh@asu.edu }

\date{}

\maketitle

\begin{abstract}
The existence, multiplicity and nonexistence of nontrivial radial convex solutions  of a system of two weakly coupled Monge-Amp\`ere equations
are established with asymptotic assumptions for an appropriately chosen parameter. The proof of the results
is based on  Krasnoselskii's fixed point theorem in a cone.
\end{abstract}


\section{Introduction}

\quad\quad In this paper we consider the existence of convex solutions to the Dirichlet problem for the weakly coupled system
\begin{equation}\label{eq3}
\begin{split}
\Big(\big(u_1'(r)\big)^N\Big)' &= \lambda N t^{N-1} f(-u_2(r)) \;\;\; \text{in} \;\; 0 < r < 1,\\
\Big(\big(u_2'(r)\big)^N\Big)' &= \lambda N t^{N-1} g(-u_1(r)) \;\;\; \text{in} \;\; 0 < r < 1,\\
 u_1'(0)=u_2'(0) & =0, \;\; u_1(1)=u_2(1)=0,
\end{split}
\end{equation}
where $N \geq 1$. A nontrivial convex solution of (\ref{eq1}) is negative on [0,1).  Such a problem arises in the study of the existence of convex radial solutions
to the Dirichlet problem for the system of the Monge-Amp\`ere equations
\begin{equation}\label{eq1}
\left\{ \begin{array}{lllll}
&\text{det} (D^2u_1)  =  \lambda f(-u_2)\;\;\textrm{in}\;\; B, \\
&\text{det} (D^2u_2)  =  \lambda g(-u_1) \;\;\textrm{in}\;\; B,\\
&u_1=u_2 =0 \;\;\textrm{on} \;\;\partial B,\;\;
\end{array} \right.
\end{equation}
where $B=\{x \in \mathbb{R}^N: |x| < 1 \}$ and $\text{det} D^2u_i$
is the determinant of  the Hessian matrix ($\frac{\partial^2u_i}{\partial x_m \partial x_n}$) of $u_i.$
For radial solution $u_i(r)$ with $r=\sqrt{\sum_1^N x_i^2}$, the Monge-Amp\`ere operator simply becomes
\begin{equation}\label{detOP}
\text{det} (D^2u_i)=\frac{(u_i')^{N-1}u_i''}{r^{N-1}}=\frac{1}{Nr^{N-1}}((u_i')^{N})'
\end{equation}
and then (\ref{eq3}) can be easily transformed into (\ref{eq1}). (\ref{detOP}) is frequently used in the literature, see
the references \cite{DELANOE1985,KUTEV1988} on radial solutions of the Monge-Amp\`ere equations  and others, e.g., Caffarelli and Li \cite{CAFFARELLILI2003}.
It can be derived from the fact that the Monge-Amp\`ere operator is rotationally invariant, see,
for example, \cite[Appendix A.2]{Goncalves2005}.

Much attention has been focused on the study of single Monge-Amp\`ere equations. Radial solutions of the boundary value problem with
a single Monge-Amp\`ere equations satisfy
\begin{equation}\label{eq-n=1}
\begin{split}
\Big(\big(u_1'(r)\big)^N\Big)'  &= \l f(-u(r))  \;\; r \in (0,1) \\
u'(0)=0, & \;\;u(1)=0.
\end{split}
\end{equation}
Kutev \cite{KUTEV1988} investigated the existence of strictly convex radial solutions of (\ref{eq-n=1}) with $f(-u)=(-u)^p$.
Delanoë \cite{DELANOE1985} treated  the existence of convex radial solutions of (\ref{eq-n=1})
for a class of more general functions, namely $\l \exp f(|x|,u,|\nabla u|)$.

By using Krasnoselskii's fixed point theorem, the author \cite{Wang2006},  and Hu and the author \cite{HU2006} showed that the existence, multiplicity and nonexistence of convex radial solutions
of (\ref{eq-n=1}) can be determined by the asymptotic behaviors of the quotient $\frac{f(u)}{u^N}$ at zero and infinity.  In a recent paper,
the author \cite{Wang2010} proved analogous results for general systems. Note that the general system in \cite{Wang2010} do not include
(\ref{eq3}).

In this paper we shall continue to establish the existence, multiplicity and nonexistence of convex radial solutions of the weakly coupled  system (\ref{eq3})
(and (\ref{eq1})) for various combinations of asymptotic behavior of $f,g$ at zero and infinity based on Krasnoselskii's fixed point theorem.
The author \cite{Wang2009} showed the existence of convex redial solutions of (\ref{eq3})
in superlinear and sublinear cases.

First,  let
$$f_0 =\lim_{x \to 0^+} \frac{f(x)}{x^N}, \quad  f_{\infty} =\lim_{x \to \infty} \frac{f(x)}{x^N},$$
and
$$g_0 =\lim_{x \to 0^+} \frac{g(x)}{x^N}, \quad  g_{\infty} =\lim_{x \to \infty} \frac{g(x)}{x^N}.$$

The main results in \cite{Wang2009} is Theorem \ref{th1}.
We include and prove it here only for completeness.

\sthm \cite{Wang2009}\mlabel{th1} Assume $f,g: [0,\infty) \to [0,\infty)$ are continuous. \\
(a). If $f_0=g_0 =0$ and $f_{\infty}=g_{\infty}=\infty$, then (\ref{eq3}) (and (\ref{eq1})) has
a nontrivial convex  solution for all $\l>0$.\\
(b). If $f_0=g_0 =\infty$ and $f_{\infty}=g_{\infty}=0$, then (\ref{eq3}) (and (\ref{eq1})) has
a nontrivial convex solution for all $\l>0$.
\ethm
Our main results in this paper are:
\sthm\mlabel{th2} Assume $f,g: [0,\infty) \to [0,\infty)$ are continuous and $f(x)>0, g(x)>0$ for $x>0$.\\
(a). If either $f_0=g_0 =0$, or $f_{\infty}=g_{\infty}=0$, then there exists a $\l_0 > 0$ such that
for all $ \l > \l_0$ (\ref{eq3}) (and (\ref{eq1}))  has a nontrivial convex solution. \\
(b). If either $f_{0}=g_{0}=\infty$, or $f_{\infty}=g_{\infty}=\infty$, then there exists a $\l_0 > 0$ such that
for all $ 0< \l < \l_0$ (\ref{eq3}) (and (\ref{eq1})) has a nontrivial convex solution. \\
(c). If $f_0=g_0=f_{\infty}=g_{\infty}=0$, then there exists a $\l_0 > 0$ such that
for all $\l > \l_0$ (\ref{eq3}) (and (\ref{eq1})) has two nontrivial convex solutions. \\
(d). If $f_0=g_0 =f_{\infty}=g_{\infty}=\infty$, then there exists a $\l_0 > 0$ such that
for all $0< \l < \l_0$ (\ref{eq3}) (and (\ref{eq1})) has two nontrivial convex solutions.\\
(e). If $f_0, g_0, f_{\infty}, g_{\infty} < \infty$, then there exists a $\l_0 > 0$ such that
for all $0< \l < \l_0$ (\ref{eq3}) (and (\ref{eq1})) has no nontrivial convex radial solution. \\
(f). If $f_0, g_0, f_{\infty}, g_{\infty}> 0$, then there exists a $\l_0 > 0$ such that
for all $\l > \l_0$ (\ref{eq3}) (and (\ref{eq1})) has no nontrivial convex radial solution.\\
\ethm

\srmark\label{rem2}
A nontrivial solution $(u_1(r),u_2(r))$ of (\ref{eq3}) has at least one nonzero component. It is possible that one of its components is zero.
$(v_i(r),v_2(r))=(-u_1(r), -u_2(r))$ is a nontrivial solution to (\ref{eq4}) or a positive fixed point of $\vect{T}_{\l}$ in (\ref{T_def}).
From the integral expression  (\ref{T_def}), we have for $r \in (0,1)$
$$v'_1 (r)=-\big(\lambda \int_0^r N \t^{N-1}f(v_2(\\t))d\t\big)^{\frac{1}{N}}$$
and
$$v'_2 (r)=-\big(\lambda \int_0^r N \t^{N-1}g(v_1(\\t))d\t\big)^{\frac{1}{N}}$$
Therefore, each component $u'_i(r)=-v'_i(r)$ is nondecreasing and $u_i=-v_i(r), i=1,2$ are convex.
\ermark

\srmark\label{rem2}
Apparently the intervals of the parameter $\lambda$ for ensuring the existence
of convex solutions of (\ref{eq3}) are not necessarily optimal. The estimates of the operator in Section \ref{preli} can be improved. We will address them
in the future.  \ermark

\section{Preliminaries}\label{preli}
With a simple transformation $v_i=-u_i, i=1,2$ (\ref{eq3}) can be brought to the following equation

\begin{equation}\label{eq4}
\left\{ \begin{array}{llll}
\Big(\big(-v_1'(r)\big)^N\Big)' &=  \lambda N r^{N-1} f(v_2),\;\; 0<r<1, \\
\Big(\big(-v_2'(r)\big)^N\Big)' &=  \lambda N r^{N-1} g(v_1),\;\; 0 < r <1,\\
 v_i'(0)=v_i(1) &= 0, \;\; i=1,2.
\end{array} \right.
\end{equation}
Now we treat positive concave classical solutions of (\ref{eq4}).

We recall some concepts and conclusions of an operator in a cone. Let $X$
be a Banach space and $K$ be a closed, nonempty subset of $X$. $K$
is said to be a cone if $(i)$~$\alpha u+\beta v\in K $ for all
$u,v\in K$ and all $\alpha,\beta \geq 0$ and $(ii)$~$u,-u\in K$ imply
$u=0$. We shall employ Krasnoselskii's fixed
point theorem to prove Theorems \ref{th1}, \ref{th2}.

\slm\mlabel{lm1} {\rm (\cite{GUOL, KRAS})} Let $X$ be a Banach
space and $K\ (\subset X)$ be a cone. Assume that $\Omega_1,\
\Omega_2$ are bounded open subsets of $X$ with $0 \in \Omega_1,\bar\Omega_1 \subset \Omega_2$, and let
$$
T: K \cap (\bar{\Omega}_2\setminus \Omega_1 ) \rightarrow K
$$
be completely continuous such that either
\begin{itemize}
\item[{\rm (i)}] $\| Tu \| \geq \| u \|,\ u\in K\cap \partial
     \Omega_1$ and $ \| Tu \| \leq \| u \|,\ u\in K\cap \partial
     \Omega_2$; or

\item[{\rm (ii)}] $\| Tu \| \leq \| u \|,\ u\in K\cap \partial
     \Omega_1$ and $\| Tu \| \geq \| u \|,\ u\in K\cap \partial
     \Omega_2$.
\end{itemize}
Then $T$ has a fixed point in $K \cap ( \bar \Omega_2 \backslash
     \Omega_1)$.

\elm

In order to apply Lemma \ref{lm1} to (\ref{eq4}), let $X$ be the Banach space
$C[0,1] \times  C[0,1]$ and, for
$(v_1,v_2) \in X,$
$$\displaystyle{\norm{(v_1,v_2)}= \norm{v_1}+ \norm{v_2}}$$
where $\norm{v_i}=\sup_{t\in[0,1]} \abs{v_i(t)}.$
Define $K$ to be a cone in $X$ by
\begin{equation*}
K = \{(v_1,v_2) \in X: v_i(t) \geq 0,\; t \in [0,1], \min\limits_{\frac{1}{4}\leq t \leq \frac{3}{4}}v_i(t) \geq \frac{1}{4} \norm{v_i}, i=1,2 \}.
\end{equation*}
For $r>0$ let
$$
\O_r  = \{(v_1,v_2) \in K: \norm{(v_1,v_2)} < r \}.
$$
Note that $\partial \O_r = \{(v_1,v_2) \in K: \norm{(v_1,v_2)}=r\}$. Further let $\vect{T}_{\lambda}: K \to X$ be a map with components $(T_{\lambda}^1, T_{\lambda}^2)$, which are defined by
\begin{equation}\label{T_def}
\begin{split}
 T_{\lambda}^1(v_1,v_2)(r) &= \int^1_r\Big( \lambda \int^s_0 N  \t^{N-1}f(v_2(\t))d\t\Big)^{\frac{1}{N}}ds, \; r \in [0,1],\\
 T_{\lambda}^2(v_1,v_2)(r) &= \int^1_r \Big( \lambda \int^s_0 N  \t^{N-1}g(v_1(\t))d\t\Big)^{\frac{1}{N}}ds, \; r \in [0,1].\\
\end{split}
\end{equation}
It is straightforward to verify that (\ref{eq4}) is equivalent to the fixed point equation
$$
\vect{T}_{\lambda}(v_1,v_2)=(v_1,v_2) \quad \text{in} \quad K.
$$
Thus, if $(v_1,v_2) \in K$ is a positive fixed point of $\vect{T}_{\lambda}$, then $(-v_1, -v_2)$ is a  convex solution of (\ref{eq3}). Conversely, if $(u_1, u_2)$ is a convex solution of (\ref{eq3}), then $(-u_1, -u_2)$ is a fixed point of $\vect{T}_{\lambda}$ in $K$.

The following lemma is a standard result due to the concavity of $v$, see e.g. \cite{Wang2003}. We prove it here only for completeness.
\slm\mlabel{lm2}
Let $v(t) \in C^1[0,1]$ for $ t \in [0,1]$. If $v(t) \geq 0 $ and $v'(t)$ is nonincreasing on $[0,1]$.
Then
$$
v(t) \geq \min\{t, 1-t\}||v||, \quad t \in [0, 1]
$$
where $||v||=\max_{t \in [0,1]}v(t).$ In particular, $$
\min_{\frac{1}{4} \leq t \leq \frac{3}{4}}v(t) \geq \frac{1}{4} ||v||.
$$
and if $v(0)=||v||$, then
$$
v(t) \geq (1-t)||v||, t \in [0,1].
$$

\elm
\pf Since $v'(t)$ is nonincreasing, we have for \mbox{$0 \leq t_0 < t < t_1 \leq 1,$}
$$
v(t)-v(t_0)=\int^t_{t_0}v'(s)ds \geq (t- t_0)v'(t)
$$
and
$$
v(t_1)-v(t)=\int^{t_1}_{t} v'(s)ds \leq (t_1 - t)v'(t),
$$
from which, we have
$$
v(t) \geq  \frac{(t_1-t)v(t_0) + (t-t_0)v(t_1)}{t_1-t_0}.
$$
Choosing $\s \in [0,1]$ such that $v(\s)=||v|| $ and considering $[t_0, t_1]$ as either of  $[0, \s]$ and $[\s, 1]$, we have
$$
v(t)  \geq  t||v|| \quad {\rm for} \quad t \in [0, \s],\\[.2cm]
$$
and
$$
v(t)  \geq (1-t)||v|| \quad {\rm for} \quad t \in [\s,1].
$$
Hence,
$$
v(t) \geq  \min\{t, 1-t\}||v||, \quad t \in [0, 1].
$$
\epf

$\vect{T}_{\lambda}(K) \subset K$ in Lemma \ref{lm-compact} is a result of Lemma \ref{lm2}. The continuity and compactness of $\vect{T}_{\lambda}$
in Lemma \ref{lm-compact}  can be verified by the standard procedures.
\slm\mlabel{lm-compact}  Assume $f,g: [0,\infty) \to [0,\infty)$ are continuous. Then $\vect{T}_{\lambda}(K) \subset K$ and $\vect{T}_{\lambda}: K \to K$ is a compact operator and continuous.
\elm

Let $$ \G = \frac{1}{4}\int^{\frac{3}{4}}_{\frac{1}{4}}\Big(\int^{s}_{\frac{1}{4}} N \t^{N-1}d\t\Big)^{\frac{1}{N}}ds> 0.$$

\slm\mlabel{f_estimate_>*} Assume $f,g: [0,\infty) \to [0,\infty)$ are continuous. Let $ (v_1,v_2) \in K $ and $\eta > 0$. If
$$f(v_2(t)) \geq (\eta v_2(t))^N \for  t \in [\frac{1}{4}, \frac{3}{4}], $$
or
$$g(v_1(t)) \geq (\eta v_1(t))^N \for  t \in [\frac{1}{4}, \frac{3}{4}], $$
then
$$
\norm{\vect{T}_{\lambda}(v_1,v_2)} \geq \lambda^{\frac{1}{N}} \G \eta \norm{v_2},
$$
or
$$
\norm{\vect{T}_{\lambda}(v_1,v_2)} \geq \lambda^{\frac{1}{N}}  \G \eta \norm{v_1},
$$
respectively.
\elm
\pf Note, from the definition of $\vect{T}_{\lambda}(v_1,v_2)$,  that $T_{\lambda}^i (v_1,v_2)(0)$ is the maximum value of $T_{\lambda}^i (v_1,v_2)$ on [0,1].
It follows that
 \begin{equation*}
 \begin{split}
 \norm{\vect{T}_{\lambda}(v_1,v_2)} & \geq \sup_{t\in[0,1]} \abs{T_{\lambda}^1(v_1,v_2)(t)} \\[.2cm]
 & \geq  \lambda^{\frac{1}{N}} \int^{\frac{3}{4}}_{\frac{1}{4}}\Big(\int^s_{\frac{1}{4}}  N \t^{N-1}f(v_2(\t))d\t\Big)^{\frac{1}{N}}ds \\[.2cm]
 & \geq  \lambda^{\frac{1}{N}} \int^{\frac{3}{4}}_{\frac{1}{4}}\Big(\int^{s}_{\frac{1}{4}} N \t^{N-1}(\eta v_2(\t))^N d\t\Big)^{\frac{1}{N}}ds \\
 & \geq  \lambda^{\frac{1}{N}}  \int^{\frac{3}{4}}_{\frac{1}{4}}\Big(\int^{s}_{\frac{1}{4}}  N \t^{N-1}( \frac{\eta}{4}\norm{v_2})^Nd\t\Big)^{\frac{1}{N}}ds \\
 & = \lambda^{\frac{1}{N}} \G \eta \norm{v_2}.
 \end{split}
 \end{equation*}
 Similarly,
\begin{equation*}
 \norm{\vect{T}_{\lambda}(v_1,v_2)}  \geq \sup_{t\in[0,1]} \abs{T_{\lambda}^2(v_1,v_2)(t)} \geq  \lambda^{\frac{1}{N}} \G \eta \norm{v_1}.
 \end{equation*}
\epf
We define new functions $\hat{f}(t), \hat{g}(t): [0, \infty) \to [0,\infty)$ by
$$\hat{f}(t) =\max \{f(v):  0 \leq v \leq t \}, \;\; \hat{g}(t) =\max \{g(v):  0 \leq v \leq t \}.$$
Note that $\hat{f}_{0}=\lim_{t \to 0} \frac{\hat{f}(t)}{t^N},$ $\hat{f}_{\infty}=\lim_{t \to \infty} \frac{\hat{f}(t)}{t^N}$
and $\hat{g}_0, \hat{g}_{\infty}$ can be defined similarly.
\slm \cite{Wang2003} \mlabel{lm6}
Assume $f,g: [0,\infty) \to [0,\infty)$ are continuous. Then $$\hat{f}_{0}=f_{0},\;\; \hat{f}_{\infty}=f_{\infty},$$
and
$$\hat{g}_{0}=g_{0}, \;\;\hat{g}_{\infty}=g_{\infty}.$$
\elm

\slm\mlabel{f_estimate_<*}
Assume $f,g: [0,\infty) \to [0,\infty)$ are continuous. Let $ r >0 $. If there exists an $\e > 0$ such that
$$
\hat{f}(r) \leq (\e r)^N, \hat{g}(r) \leq (\e r)^N,
$$
then
$$
\norm{\vect{T}_{\lambda}(v_1,v_2)} \leq 2 \e  \lambda^{\frac{1}{N}} \norm{(v_1,v_2)} \;\; {\rm for} \;\; (v_1,v_2) \in \partial\O_{r}.
$$
\elm
\pf From the definition of $T$, for $(v_1,v_2) \in \partial\O_{r}$, we have
 \begin{eqnarray*}
 \norm{\vect{T}_{\lambda}(v_1,v_2)} & =&  \sum_{i=1}^2 \sup_{t\in[0,1]} \abs{T_{\lambda}^i(v_1,v_2)(t)}\\
  & \leq & \lambda^{\frac{1}{N}} (\int^1_0 N \t^{N-1} f(v_2(\t))d\t)^{\frac{1}{N}}+\lambda^{\frac{1}{N}} (\int^1_0 N \t^{N-1} g(v_1(\t))d\t)^{\frac{1}{N}} \\
  & \leq & \lambda^{\frac{1}{N}}(\int^1_0 N \t^{N-1} \hat{f}(r)d\t)^{\frac{1}{N}}+\lambda^{\frac{1}{N}}(\int^1_0 N \t^{N-1} \hat{g}(r)d\t)^{\frac{1}{N}} \\
  & \leq & \lambda^{\frac{1}{N}} (\int^1_0 N \t^{N-1} d\t)^{\frac{1}{N}}\e r +\lambda^{\frac{1}{N}}(\int^1_0 N \t^{N-1} d\t)^{\frac{1}{N}}\e r \\
  & = &  2 \lambda^{\frac{1}{N}} \e \norm{(v_1,v_2)}.
  \end{eqnarray*}
\epf

The following two lemmas are weak forms of Lemmas \ref{f_estimate_>*} and \ref{f_estimate_<*}.
\slm\mlabel{f_estimate_>*_weak} Assume $f,g: [0,\infty) \to [0,\infty)$ are continuous and $f(u),g(u)>0$ for $u>0$.
Let $r>0, (v_1,v_2) \in \O_r $. Then
$$
\norm{\vect{T}_{\lambda}(v_1,v_2)} \geq  4 \l ^{\frac{1}{N}} \G (\hat{m}_r)^{\frac{1}{N}},
$$
where $\hat{m}_r=\min_{\frac{r}{8} \leq t \leq r}\{ \min\{f(t), g(t)\}\}>0.$
\elm
\pf  If $ (v_1,v_2) \in \partial \O_{r}$, then $\sup_{t \in [0,1]} v_i \geq \frac{1}{2}r$ must hold for either $i=1$ or $i=2.$
If $\sup_{t \in [0,1]} v_2 \geq \frac{1}{2}r,$ then $$ \min_{\frac{1}{4}\leq t \leq \frac{3}{4}} v_2(t) \geq \frac{1}{4}\sup_{t \in [0,1]} v_2 \geq \frac{1}{8} r$$
which implies that
$$f(v_2(t)) \geq \hat{m}_r \; \rm{for} \; t \in [\frac{1}{4}, \frac{3}{4}].
$$
Otherwise, we have
$$g(v_1(t)) \geq \hat{m}_r \; \rm{for} \; t \in [\frac{1}{4}, \frac{3}{4}].
$$
It is easy to see that this lemma can be shown in a similar manner as in Lemma \ref{f_estimate_>*}.
\epf
\slm\mlabel{f_estimate_<*_weak} Assume $f,g: [0,\infty) \to [0,\infty)$ are continuous and $f(u),g(u)>0$ for $u>0$.
Let $r>0, (v_1,v_2) \in \O_r $. Then
$$
\norm{\vect{T}_{\lambda}(v_1,v_2)} \leq 2 \l ^{\frac{1}{N}} (\hat{M}_r)^{\frac{1}{N}},
$$
where $\hat{M}_r=\max\{ f(t)+g(t):  0 \leq t \leq r \}>0$
\elm
\pf Since $ f(v_1(t)), f(v_1(t)) \leq \hat{M}_r = \rm{for}\; t \in [0,1]$, it is easy to see that this lemma can be shown in a similar manner as in Lemma \ref{f_estimate_<*}.
\epf

\section{Proof of Theorem \ref{th1}}
\pf
Part (a). It follows from Lemma ~\ref{lm6} that $\hat{f}_0=0, \hat{g}_0=0.$ Therefore, we can choose $r_1 > 0$
so that $\hat{f}(r_1) \le (\e r_1)^N, \hat{g}(r_1) \le (\e r_1)^N$ where the constant $\e> 0$ satisfies
$$
2 \e  \lambda^{\frac{1}{N}} <1 .
$$
We have by Lemma ~\ref{f_estimate_<*} that
$$
\norm{\vect{T}_{\lambda}(v_1,v_2)} \leq  2  \lambda^{\frac{1}{N}} \e \norm{(v_1,v_2)} < \norm{(v_1,v_2)} \quad \textrm{for} \quad  (v_1,v_2) \in \partial\O_{r_1}.
$$
Now, since $f_{\infty} = \infty, g_{\infty} = \infty$, there is
an $\hat{H} > 0$ such that
$$
f(v) \geq (\eta v)^{N}, g(v) \geq (\eta v)^{N}
$$
for $ v \geq \hat{H}$ ,
where $\eta > 0$ is chosen so that
$$
\frac{1}{2} \lambda^{\frac{1}{N}} \G \eta > 1.
$$
Let $r_2 = \max\{2r_1,8\hat{H} \}$. If $ (v_1,v_2) \in \partial \O_{r_2}$, $\sup_{t \in [0,1]} v_i \geq \frac{1}{2}r_2$
must hold for either $i=1$ or $i=2$. Without loss of generality, assume that
$\sup_{t \in [0,1]} v_1 \geq \frac{1}{2}r_2.$
Then $$ \min_{\frac{1}{4}\leq t \leq \frac{3}{4}} v_1(t) \geq \frac{1}{4}\sup_{t \in [0,1]} v_1 \geq \frac{1}{8} r_2\geq \hat{H},$$
which implies that
$$g(v_1(t)) \geq (\eta v_1(t))^N \; \rm{for} \; t \in [\frac{1}{4}, \frac{3}{4}].
$$
It follows from Lemma ~\ref{f_estimate_>*} that
$$
\norm{\vect{T}_{\lambda}(v_1,v_2)} \geq \G  \lambda^{\frac{1}{N}} \eta \norm{v_1} \geq \frac{1}{2}  \lambda^{\frac{1}{N}} \G \eta r_2 > r_2 =\norm{(v_1,v_2)}.
$$
By Lemma ~\ref{lm1}, $\vect{T}_{\lambda}$ has a fixed point $(v_1,v_2) \in  \O_{r_2} \setminus \bar{\O}_{r_1}$.
The fixed point $(v_1,v_2) \in \O_{r_2} \setminus \bar{\O}_{r_1}$ is the desired positive solution of (\ref{eq4}).

Part (b). Since $f_{0} = \infty, g_{0} = \infty$, there is
an $r_1 > 0$ such that
$$
f(v) \geq (\eta v)^{N}, g(v) \geq (\eta v)^{N}
$$
for $0 <  v \leq  r_1$ ,
where $\eta > 0$ is chosen so that
$$
\lambda^{\frac{1}{N}} \G \eta > 1.
$$
If $ (v_1,v_2) \in  \partial \O_{r_1}$, then
$$f(v_2(t)) \geq (\eta v_2)^N, \;\; g(v_1(t)) \geq (\eta v_1)^N {\rm for } \;\; t \in [0,1].$$
Lemma ~\ref{f_estimate_>*} implies that
$$
\norm{\vect{T}_{\lambda}(v_1,v_2)} \geq \lambda^{\frac{1}{N}} \G \eta \norm{(v_1,v_2)} > \norm{(v_1,v_2)} \quad \textrm{for}\quad  (v_1,v_2) \in \partial\O_{r_1}.
$$
We now determine $\O_{r_2}$. It follows from Lemma ~\ref{lm6}
that $\hat{f}_{\infty}=0$ and $\hat{g}_{\infty}=0.$
Therefore there is an $r_2>2r_1$ such that
$$
\hat{f}(r_2) \le (\e r_2)^N,\;\; \hat{g}(r_2) \le (\e r_2)^N ,
$$
where the constant $ 2 \lambda^{\frac{1}{N}} \e <1.$
Thus, we have by Lemma ~\ref{f_estimate_<*} that
$$
\norm{\vect{T}_{\lambda}(v_1,v_2)} \leq 2 \e \lambda^{\frac{1}{N}} \norm{(v_1,v_2)} < \norm{(v_1,v_2)} \quad \textrm{for}\quad  (v_1,v_2) \in \partial\O_{r_2}.
$$
By Lemma ~\ref{lm1}, $\vect{T}_{\lambda}$ has a fixed point $(v_1, v_2)$ in  $\O_{r_2} \setminus \bar{\O}_{r_1}$.
And $(v_1, v_2)$ is the desired positive solution of (\ref{eq4}).
\epf

\section{Proof of Theorem \ref{th2}}
\pf
Part (a).
Fix a number $r_1 > 0$. Lemma \ref{f_estimate_>*_weak} implies that there exists a $\l_0 >0$ such that
$$
\norm{\vect{T}_{\lambda}(v_1,v_2)} > \norm{(v_1,v_2)}=r_1, \;\; \rm{for}\;\; (v_1,v_2) \in  \partial \O_{r_1}, \l > \l_0.
$$
If $f_0=g_0 =0$, it follows from Lemma \ref{lm6} that
$$
\hat{f}_0=\hat{g}_0=0.
$$
Therefore, we can choose $0 < r_2 < r_1 $ so that
$$
\hat{f}(r_2) \le (\e r_2)^N, \hat{g}(r_2) \le (\e r_2)^N
$$
where the constant $\e> 0$ satisfies
$$
2 \e  \lambda^{\frac{1}{N}} <1.
$$
We have by Lemma ~\ref{f_estimate_<*} that
$$
\norm{\vect{T}_{\lambda}(v_1,v_2)} \leq  2  \lambda^{\frac{1}{N}} \e \norm{(v_1,v_2)} < \norm{(v_1,v_2)} \quad \textrm{for} \quad  (v_1,v_2) \in \partial\O_{r_2}.
$$
If $f_{\infty}=g_{\infty}=0$, it follows from Lemma \ref{lm6} that then $\hat{f}_{\infty}=\hat{g}_{\infty}=0$.
Therefore there is an $r_3>2r_1$ such that
$$
\hat{f}(r_3) \le (\e r_3)^N, \hat{g}(r_3) \le (\e r_3)^N
$$
where the constant $\e > 0$ satisfies
$$
2 \e  \lambda^{\frac{1}{N}} <1.
$$
Thus, we have by Lemma ~\ref{f_estimate_<*} that
$$
\norm{\vect{T}_{\lambda}(v_1,v_2)} \leq  2  \lambda^{\frac{1}{N}} \e \norm{(v_1,v_2)} < \norm{(v_1,v_2)} \quad \textrm{for} \quad  (v_1,v_2) \in \partial\O_{r_3}.
$$
It follows from Lemma ~\ref{lm1} that $\vect{T}_{\l}$ has a fixed point in  $\O_{r_1} \setminus \bar{\O}_{r_2}$ or $\O_{r_3} \setminus \bar{\O}_{r_1}$.
Consequently, (\ref{eq4}) has a positive solution for $ \l > \l_0$.

Part (b).
Fix a number $r_1 > 0$. Lemma \ref{f_estimate_<*_weak} implies that there exists a $\l_0 >0$ such that
$$
\norm{\vect{T}_{\lambda}(v_1,v_2)} < \norm{(v_1,v_2)}=r_1, \; {\rm for} \; \vect{v} \in  \partial \O_{r_1},\; 0< \l < \l_0.
$$
If $f_{0}=g_{0} = \infty$, there is an $0<r_2<r_1$ such that
$$
f(v) \geq (\eta v)^{N}, g(v) \geq (\eta v)^{N}
$$
for $0 <  v \leq  r_2$ ,
where $\eta > 0$ is chosen so that
$$
\lambda^{\frac{1}{N}} \G \eta > 1.
$$
If $ (v_1,v_2) \in  \partial \O_{r_2}$, then
$$f(v_2(t)) \geq (\eta v_2)^N, \;\; g(v_1(t)) \geq (\eta v_1)^N {\rm for } \;\; t \in [0,1].$$
Lemma ~\ref{f_estimate_>*} implies that
$$
\norm{\vect{T}_{\lambda}(v_1,v_2)} \geq \lambda^{\frac{1}{N}} \G \eta \norm{(v_1,v_2)} > \norm{(v_1,v_2)} \quad \textrm{for}\quad  (v_1,v_2) \in \partial\O_{r_2}.
$$

If $f_{\infty}=g_{\infty} = \infty$, there is
an $\hat{H} > 0$ such that
$$
f(v) \geq (\eta v)^{N}, g(v) \geq (\eta v)^{N}
$$
for $ v \geq \hat{H}$ ,
where $\eta > 0$ is chosen so that
$$
\frac{1}{2} \lambda^{\frac{1}{N}} \G \eta > 1.
$$
Let $r_3 = \max\{2r_1,8\hat{H} \}$. If $ (v_1,v_2) \in \partial \O_{r_3}$,  $\sup_{t \in [0,1]} v_i \geq \frac{1}{2}r_3$ must hold for
either $i=1$ or $i=2$. Without loss of generality, assume that
$\sup_{t \in [0,1]} v_1 \geq \frac{1}{2}r_3.$
Then $$ \min_{\frac{1}{4}\leq t \leq \frac{3}{4}} v_1(t) \geq \frac{1}{4}\sup_{t \in [0,1]} v_1 \geq \frac{1}{8} r_3\geq \hat{H},$$
which implies that
$$g(v_1(t)) \geq (\eta v_1(t))^N \; \rm{for} \; t \in [\frac{1}{4}, \frac{3}{4}].
$$
It follows from Lemma ~\ref{f_estimate_>*} that
$$
\norm{\vect{T}_{\lambda}(v_1,v_2)} \geq \G  \lambda^{\frac{1}{N}} \eta \norm{v_1} \geq  \frac{1}{2}  \lambda^{\frac{1}{N}} \G \eta r_3 > r_3 =\norm{(v_1,v_2)}.
$$
It follows from Lemma ~\ref{lm1} that
$\vect{T}_{\l}$ has a fixed point in  $\O_{r_1} \setminus \bar{\O}_{r_2}$ or $\O_{r_3} \setminus \bar{\O}_{r_1}.$
Consequently, (\ref{eq4}) has a positive solution for $ 0< \l < \l_0$.

Part (c).
Fix two numbers $0 < r_{3} < r_{4}.$ Lemma \ref{f_estimate_>*_weak} implies that there exists a $\l_0 >0$ such that for $\l >  \l_{0}$,
$$
\norm{\vect{T}_{\l}(v_1,v_2)} > \norm{(v_1,v_2)}, \;\;\rm{for} \;\; (v_1,v_2) \in  \partial
\O_{r_i}, \;\; (i=3,4).
$$
Since $f_0=g_0=f_{\infty}=g_{\infty}=0$, it follows from the proof of Theorem ~\ref{th2} (a) that we can choose $0< r_1 < r_3/2$ and $r_2>2r_4$  such that
$$\norm{\vect{T}_{\l}(v_1,v_2)}  <  \norm{(v_1,v_2)}, \;\; {\rm for} \;\; (v_1,v_2) \in \partial \O_{r_i}, \;\; (i=1,2).$$
It follows from Lemma ~\ref{lm1} that
 $\vect{T}_{\l}$ has two fixed points $(v_1,v_2)$ and $(u_1,u_2)$ such that  $(v_1,v_2) \in \O_{r_3} \setminus \bar{\O}_{r_1}$ and
 $(u_1,u_2) \in \O_{r_2} \setminus \bar{\O}_{r_4}$ , which are
the desired distinct positive solutions of (\ref{eq4}) for $ \l > \l_0$ satisfying
$$
r_1 < \norm{(v_1,v_2)} < r_3 < r_4 < \norm{(u_1,u_2)} < r_2.
$$

Part (d).
Fix two numbers $0 < r_{3} < r_{4}.$ Lemma \ref{f_estimate_<*_weak} implies that there exists a $\l_0 >0$ such that for $ 0< \l < \l_{0} $,
$$
\norm{\vect{T}_{\l}(v_1,v_2)} < \norm{(v_1,v_2)}, \;\;\rm{for} \;\; (v_1,v_2) \in  \partial
\O_{r_i}, \;\; (i=3,4).
$$
Since $f_0=g_0=f_\infty=g_\infty=\infty$, it follows from the proof of Theorem ~\ref{th2} (b) that
we can choose $0< r_1 < r_3/2$ and $r_2>2r_4$ such that
$$\norm{\vect{T}_{\l}(v_1,v_2)} > \norm{(v_1,v_2)}, \;\; {\rm for} \;\; (v_1,v_2) \in  \partial \O_{r_i}, \;\; (i=1,2).$$
It follows from Lemma ~\ref{lm1} that
$\vect{T}_{\l}$ has two fixed points $(v_1,v_2)$ and $(u_1,u_2)$ such that  $(v_1,v_2) \in \O_{r_3} \setminus \bar{\O}_{r_1}$
and $(u_1,u_2) \in \O_{r_2} \setminus \bar{\O}_{r_4}$ ,
which are the desired distinct positive solutions of (\ref{eq4}) for $ \l < \l_0$ satisfying
$$
r_1 < \norm{(v_1,v_2)} < r_3 < r_4 < \norm{(u_1,u_2)} < r_2.
$$

Part (e). Since $f_0,g_0, f_{\infty},g_{\infty}< \infty$, it follows that $\hat{f_0},\hat{g_0}, \hat{f}_{\infty},\hat{g}_{\infty}< \infty$
and there exist positive numbers $\e_1^i$, $ \e_2^i$, $r_1^i$ and $r_2^i$ such that $r_1^i < r_2^i, i=1,2$,
$$
  \hat{f}(v)  \leq   \e_1^1 v ^{N}\; \rm{for } \;0 \leq v \leq r_1^1, \;\;\;  \hat{g}(v)  \leq   \e_1^2 v ^{N}\; \rm{for } \; 0 \leq v \leq r_1^2,
$$
and
$$
  \hat{f}(v)  \leq   \e_2^1 v ^{N}\; \rm{for } \; v \geq r_2^1, \;\;\;  \hat{g}(v)  \leq   \e_2^2 v ^{N}\; \rm{for } \; v \geq r_2^2,
$$
Let
$$
\e^1 = \max\{\e_1^1, \e_2^1, \max\{\frac{\hat{f}(v)}{v^N}:  r_1^1 \leq v \leq  r_2^1 \}\} > 0
$$
$$
\e^2 = \max\{\e_1^2, \e_2^2, \max\{\frac{\hat{g}(v)}{v^N}:  r_1^2 \leq v \leq  r_2^2 \}\} > 0
$$
and $ \e = (\max\limits_{i=1,2} \{\e^i\})^{\frac{1}{N}} > 0.$ Thus, we have
$$
 \hat{f}(v) \leq  (\e v)^N,\;\; \hat{g}(v) \leq  (\e v)^N \;\rm{ for }\; v \geq 0.
$$
Assume $(v_1(t),v_2(t))$ is a positive solution of (\ref{eq4}). We will show that this leads to a contradiction
for $0<  \l < \l_0,$
where
$$
\l_0=(\frac{1}{2\e})^N.
$$
In fact, for $0<  \l < \l_0$, since  $\vect{T}_{\l}(v_1(t),v_2(t)) = (v_1(t),v_2(t))$ for $ t \in [0,1]$, by Lemma ~\ref{f_estimate_<*}, we have
 \begin{eqnarray*}
  \norm{(v_1,v_2)} & = &  \norm{\vect{T}_{\l}(v_1,v_2)}\\
  & \leq &  2 \e  \lambda^{\frac{1}{N}} \norm{(v_1,v_2)} \\
  & < &    \norm{(v_1,v_2)}
\end{eqnarray*}
which is a contradiction.

Part (f). Since $f_0, g_0, f_{\infty}, g_{\infty}> 0$, there exist positive numbers $\e_1^i$, $ \e_2^i$, $r_1^i$ and $r_2^i$ such that $r_1^i < r_2^i, i=1,2$,
$$
  f(v)  \geq   \e_1^1 v ^{N}\; \rm{for } \;0 \leq v \leq r_1^1, \;\;\;  g(v)  \geq   \e_1^2 v ^{N}\; \rm{for } \; 0 \leq v \leq r_1^2,
$$
and
$$
  f(v)  \geq   \e_2^1 v ^{N}\; \rm{for } \; v \geq r_2^1, \;\;\;  g(v)  \geq   \e_2^2 v ^{N}\; \rm{for } \; v \geq r_2^2,
$$
Let
$$
\e^1 = \min\{\e_1^1, \e_2^1, \min\{\frac{f(v)}{v^N}:  r_1^1 \leq v \leq  r_2^1 \}\} > 0
$$
$$
\e^2 = \min\{\e_1^2, \e_2^2, \min\{\frac{g(v)}{v^N}:  r_1^2 \leq v \leq  r_2^2 \}\} > 0
$$
and $ \e = (\min \limits_{i=1,2} \{\e^i\})^{\frac{1}{N}} > 0.$ Thus, we have
$$
 f(v) \geq  (\e v)^N,\;\; g(v) \geq  (\e v)^N \;\rm{ for }\; v \geq 0.
$$
Assume $(v_1(t),v_2(t))$ is a positive solution of (\ref{eq4}). We will show that this leads to a contradiction
for $\l > \l_0 =(\frac{1}{ \Gamma \e })^N$.  In fact, since  $\vect{T}_{\l}(v_1(t),v_2(t)) = (v_1(t),v_2(t))$ for $ t \in [0,1]$, it
follows from Lemma ~\ref{f_estimate_>*} that, for $\l > \l_0$,
 \begin{eqnarray*}
  \norm{(v_1,v_2)} & = & \norm{\vect{T}_{\l}(v_1,v_2)} \\[.2cm]
  & \geq &   \Gamma \e  \lambda^{\frac{1}{N}} \norm{(v_1,v_2)} \\[.2cm]
  & > & \norm{(v_1,v_2)},
\end{eqnarray*}
which is a contradiction.\epf

\end{document}